\documentclass[12pt]{article}

\usepackage{graphicx}
\usepackage{latexsym,amssymb}
\usepackage{amsthm}
\usepackage{indentfirst}
\usepackage{amsmath}
\usepackage{color}
\usepackage{xcolor}
\usepackage{fourier}
\usepackage[all]{xy}
\usepackage[colorlinks=true,backref=page]{hyperref}

\newtheorem{theoremalph}{Theorem}

\newtheorem*{Main Theorem}{Main Theorem}

\newtheorem{Theorem}{Theorem}[section]
\newtheorem*{Theorem A}{Theorem A}
\newtheorem*{Theorem A'}{Theorem A'}
\newtheorem*{Theorem B'}{Theorem B'}

\newtheorem{Proposition}[Theorem]{Proposition}
\newtheorem{Lemma}[Theorem]{Lemma}

\newtheorem{Question}{Question}

\newtheorem*{Notation}{Notation}
\newtheorem{Remark}[Theorem]{Remark}

\newtheorem{Corollary}[Theorem]{Corollary}

\newtheorem{Claim-numbered}[Theorem]{Claim}

  \def\cZ{{\mathbb C}}

 \def\NN{{\mathbb N}} 

\def\QQ{{\mathbb Q}} \def\RR{{\mathbb R}} \def\SS{{\mathbb S}}
\def\TT{{\mathbb T}}

 \def\ZZ{{\mathbb Z}}

  \def\cG{{\cal G}}  
    
\def\cZ{{\cal C}}    \def\cU{{\cal U}}
    
\def\cE{{\cal E}}    \def\cW{{\cal W}}
\def\cF{{\cal F}}    
  \def\cZ{{\cal Z}}

\newcommand{\Id}{\operatorname{Id}}

\newcommand{\diff}{{\operatorname{Diff}}}

\def\diff{\operatorname{Diff}}

\def\dim{\operatorname{dim}}

\def\ud{\operatorname{d}}
\def\e{{\varepsilon}}

\def\det{\operatorname{det}}

\def\homeo{\operatorname{Homeo}}
\def\Acc{\operatorname{Acc}}
\def\wh{\widehat}

\begin{document}

\title{
Centralizers of  derived-from-Anosov systems on $\TT^3$: rigidity versus triviality
}

\author{Shaobo Gan, Yi Shi, Disheng Xu and  Jinhua Zhang}


\maketitle

\begin{abstract}
In this paper, we study the centralizer of a partially hyperbolic diffeomorphism on $\TT^3$ which is homotopic to an Anosov automorphism, and we show that either its centralizer is virtually trivial or such diffeomorphism is smoothly conjugate to its linear part.

\hspace{-1cm}\mbox
\smallskip

\end{abstract}


\section{Introduction}
Let $M$ be a closed $C^{\infty}$ Riemannian manifold.
Given a  diffeomorphism $f\in\diff^r(M)$ and $r\in[1,\infty], s\in [0,r]$,   the \emph{$C^s$-centralizer} of $f$ is defined as   $$\cZ^s(f)=\big\{g\in\diff^s(M): g\circ f=f\circ g \big\}.$$ 
By definition, $f$ is $C^s$ conjugate to itself by $g\in\cZ^s(f)$ and the centralizer $\cZ^s(f)$ always contains the cyclic group $<f>:=\{f^n\}_{n\in\ZZ}$ generated by $f$. One says that $f$ has \emph{trivial $C^s$-centralizer} if $\cZ^s(f)=\{f^n\}_{n\in\ZZ}$, and $f$ has \emph{virtually trivial $C^s$-centralizer} if the cyclic group $<f>$ is a finite index subgroup of $\cZ^s(f)$.  

The centralizers of diffeomorphisms play important role in several topics of dynamical systems. For instance, people attempt to classify diffeomorphisms up to differentiable conjugacies, especially in the study of circle diffeomorphisms \cite{He}. On the other hand, the property of the centralizer of a diffeomorphism can give some consequence in  foliation theory, see for instance~\cite{B}. Moreover, the centralizer of diffeomorphisms  is closely related to the study of higher rank abelian actions on manifolds, see for instance~\cite{DK,HW}. 
Smale~\cite{Smale90,Smale98} proposed the following conjecture, and considered it as one of the mathematical problems in this century:

\medskip
\noindent{\bf Conjecture.} \emph{There exists a residual subset ${\cal U}\subseteq\diff^r(M)$, such that every $f\in{\cal U}$ has trivial $C^r$-centralizer.}
 
\medskip

This conjecture has been proved in the case $r=1$, see~\cite{BCW,BCVW}. For $r>1$, it has been proved by Palis and Yoccoz \cite{PY1,PY2} that $C^r$-open dense Anosov diffeomorphisms have trivial centralizer. See also \cite{RV} for hyperbolic basic sets.  

Since Pugh and Shub proposed the famous Stable Ergodicity Conjecture, the partially hyperbolic system has been one of the main topics of research in dynamical systems. 
A diffeomorphism    $f\in\diff^r(M)$ is \emph{partially hyperbolic} if there exist a $Df$-invariant splitting $TM=E^{s}\oplus E^c\oplus E^{u}$ and $N\in\NN$ such that 
\begin{itemize}
	\item \textbf{Uniform contracting and expanding}: for any $x\in M$, one has 
	$$\|Df^{N}|_{E^{s}(x)}\|\leq \frac{1}{2}
	\qquad \textrm{ and } \qquad
	\|Df^{-N}|_{E^{u}(x)}\|\leq \frac{1}{2};$$
	\item \textbf{Domination}: for any $x\in M$, one has 
	$$\|Df^{N}|_{E^{s}(x)}\|\cdot\|Df^{-N}|_{E^c(f^N(x))}\|\leq \frac{1}{2}
	\quad \textrm{ and } \quad
	\|Df^N|_{E^c(x)}\|\cdot \|Df^{-N}|_{E^{u}(f^N(x))}\|\leq \frac 12.$$
\end{itemize}
It is clear that the set consists of all $C^r$ partially hyperbolic diffeomorphisms is an open subset of $\diff^r(M)$ in the $C^r$-topology.
If either the strong stable bundle $E^s$, or the strong unstable bundle $E^u$ is trivial, we say $f$ is \emph{weakly partially hyperbolic}.

Not every manifold support a partially hyperbolic diffeomorphism. For instance, there are no partially hyperbolic diffeomorphisms on $\SS^3$ \cite{BI}. The simpliest manifold supporting a partially hyperbolic diffeomorphism is the 3-torus $\TT^3$. 
It has been proved \cite{BBI2,Po} that if $f\in\diff^r(\TT^3)$ is partially hyperbolic,
then its linear part $L_f:\pi_1(\mathbb{T}^3)=\mathbb{Z}^3\rightarrow\mathbb{Z}^3$ induced by $f$ on the fundamental group of $\TT^3$ is also partially hyperbolic,
that is, $L_f\in{\rm GL}(3,\mathbb{Z})$ has three real eigenvalues with different moduli: one has modulus larger than $1$ and one has modulus smaller than one.
So there are two classes of partially hyperbolic diffeomorphisms on $\mathbb{T}^3$:
\begin{itemize}
	\item either $L_f\in{\rm GL}(3,\mathbb{Z})$ has an eigenvalue with modulus 1;
	\item or $L_f\in{\rm GL}(3,\mathbb{Z})$ is Anosov.
\end{itemize}
Since $L_f\in{\rm GL}(3,\mathbb{Z})$ defines a diffeomorphism on $\TT^3$, we still denote $L_f\in\diff^r(\TT^3)$ for simplicity. It is clear that $f$ is homotopic to $L_f$ on $\TT^3$. 

It has been proved by Potrie \cite{Po} that if $f$ is in the first class, and the non-wandering set of $f$ is whole $\TT^3$, then the center foliation of $f$ forms a circle bundle. 
We say such $f$ is a partially hyperbolic skew-product diffeomorphism on $\TT^3$. Recently, the centralizers of certain classes of conservative partially hyperbolic skew-product diffeomorphisms have been totally classified by \cite{DWX}. 

If $L_f\in{\rm GL}(3,\mathbb{Z})$ is Anosov, i.e. $f$ is homotopic to an Anosov automorphism, then we say $f$ is a \emph{derived-from-Anosov diffeomorphism}, or a \emph{DA-diffeomorphism}. In this case, $L_f$ is an Anosov automorphism with real simple spectrum.
The study of partially hyperbolic derived-from-Anosov diffeomorphisms originated from Ma\~n\'e \cite{M}. He constructed a robustly topologically mixing 
partially hyperbolic DA-diffeomo\-rphism which is not Anosov. 

The partially hyperbolic DA-diffeomorphisms have been studied extensively, in both topological and measure-theoretic  aspects. For instance, every DA-diffeomorp\-hism $f$ is dynamically coherent and leaf conjugate to its linear part $L_f$ \cite{BBI1,Ha,HaPo0,Po}. It has been shown in \cite{BFSV, U} that partially hyperbolic DA-diffeomorphisms are intrinsically ergodic. The disintegration of an invariant measure along the center foliation has been studied
in detail \cite{VY}. Moreover, every conservative partially hyperbolic DA-diffeomorphism is ergodic \cite{HaU,GS}. In certain settings, it is further known to be Bernoulli \cite{PTV}.

\medskip

In this  paper, we study the centralizer of these diffeomorphisms.

  \begin{theoremalph}~\label{thm.main-thm}
   Let $f\in\diff^\infty(\TT^3)$ be a partially hyperbolic derived-from-Anosov diffeomorphism.
   Then one has the following dichotomy:
  	\begin{itemize}
  		\item either the   $C^\infty$ centralizer of $f$ is virtually trivial, and 
  		$$
  		\#\big\{g\in\cZ^\infty(f): \textrm{$g$ is homotopic to identity } \big\}\leq |\det(L_f-\Id)|;
  		$$
  		\item  or $f$ is $C^\infty$-conjugate to $L_f$, thus $\cZ^\infty(f)\cong\cZ^\infty(L_f)$. 
  	\end{itemize}
  \end{theoremalph}

For the centralizer of the partially hyperbolic derived-from-Anosov diffeomorphisms with lower regularity, we get the following theorem. And Theorem \ref{thm.main-thm} is a direct corollary of this theorem.

\begin{theoremalph}~\label{thm.main-thm-low-regularity}
Let $f\in\diff^r(\TT^3)$ ($r>1$) be a partially hyperbolic derived-from-Anosov diffeomorphism.
Then $f$ satisfies one of the following properties:
\begin{itemize}
		\item the $C^r$ centralizer of $f$ is virtually trivial and 
		$$\#\big\{g\in\cZ^r(f): \textrm{$g$ is homotopic to identity } \big\}\leq |\det(L_f-\Id)|;$$
		\item $f$ is Anosov, and $C^{r-\e}$-conjugate to $L_f$ for any $\e>0$.
\end{itemize}
\end{theoremalph}

\begin{Remark}
\begin{enumerate}
	\item In the second case of Theorem \ref{thm.main-thm-low-regularity}, the loss of regularity comes from Journ\'e's theorem \cite{J}. If $r$ is not an integer, then $f$ is $C^r$-conjugate to $L_f$. Thus $\cZ^r(f)\cong\cZ^r(L_f)=\cZ^\infty(L_f)$, which is virtually $\ZZ^2$. 
	\item  	From  the proof, we can see that the dichotomy in both theorems is whether $f$ is accessible. If $f$ is accessible, then $\cZ^r(f)$ is virtually trivial. Otherwise, $f$ is forced to be smoothly  conjugate to its linear part $L_f$.
	\item 	The reason we discuss virtual triviality rather than triviality of the centralizer is that this property is more likely to be robust, for example the class of systems which satisfies the dichotomy in our paper forms an open subset in the group of $C^r$ diffeomorphisms on $\TT^3$. The readers can refer to \cite{DWX} and references therein   for more  results  on ``virtual triviality of centralizer or rigidity" for partially hyperbolic diffeomorphisms.
\end{enumerate}

\end{Remark}


A direct corollary is $C^r$-open densely, the $C^r$-centralizer of a partially hyperbolic diffeomorphism which is homotopic to an Anosov automorphism on $\TT^3$ is virtually trivial.
In particular, our result implies that the centralizer of every diffeomorphism constructed by Ma\~n\'e in \cite{M} is virtually trivial.

We don't know the case for $r=1$. Our proof strongly relies on the recent results in  \cite{GS,HaS} where $r>1$ is crucial in their arguments.
\begin{Question}
	Does Theorem~\ref{thm.main-thm-low-regularity} hold when $r=1$?
\end{Question}
\begin{Remark}
If $f$ is $C^1$-smooth and accessible, then $\cZ^1(f)$ is virtually trivial. See Remark~\ref{r.homotopic-to-id-C1-case} and Theorem~\ref{thm.joint-integrability} in Section~\ref{s.proof-of-thms}.
\end{Remark}

Recently, Barthelm\'e and Gogolev \cite{BG} studied the centralizer of volume preserving partially hyperbolic diffeomorphisms which are homotopic to identity on Seifert fibered and hyperbolic 3-manifolds. They proved that  for these diffeomorphisms, either their centralizers are virtually trivial, or they can be embedded into a $C^1$ flow.  In~\cite{Bur},  the author showed that among the set of $C^r$ ($r\geq 1$) partially hyperbolic diffeomorphisms on a closed manifold, there exists a $C^1$-open and $C^1$-dense subset where each element has discrete centralizer. 
We propose the following question.

\begin{Question}
Suppose that $M$ is a closed  3-manifold and ${\rm PH}^r(M)$ is the set of $C^r$ partially hyperbolic diffeomorphisms on $M$. Let $\cU^r\subset {\rm PH}^r(M)$ be defined as the following: $$\cU^r:=\big\{f\in {\rm PH}^r(M): \cZ^r(f) \text{~is trivial~}\big\}.$$
Does $\cU^r$  contain  an open dense subset in ${\rm PH}^r(M)$ for every $r>1$?
\end{Question}

In this paper, the linear part of the partially hyperbolic diffeomorphisms have real simple spectra. It would be interesting to know the case where the weakly partially hyperbolic diffeomorphisms  homotopic to a linear Anosov on $\TT^3$ exhibiting  complex eigenvalues.


\medskip

\noindent{\bf Acknowledgments: } We would like to thank A. Gogolev for useful comments. D.X. would like to thank Professor Amie Wilkinson for many useful discussions in the study of the general centralizer problems, and  Claim \ref{Cl: Amie} is partially inspired by discussion with her. Part of this work is done when D.X is served as Marie Curie fellow in Imperial College London. S.G. is supported by NSFC 11771025 and 11831001.  Y.S. is supported by NSFC 11701015 and 11831001. J.Z. is supported by starting grant from Beihang University. 

 \section{Preliminaries}
 
 In this section, we collect the notions and results used in this paper.
 \subsection{Dominated splitting}

 A $Df$-invariant splitting $TM=E\oplus F$ is \emph{dominated}, if there exist  $C>1$ and $\lambda\in(0,1)$ such that 
 $$\|Df^n|_{E(x)}\|\cdot \|Df^{-n}|_{F(f^n(x))}\|\leq  C\cdot \lambda^n \textrm { for any $
 x\in M$ and any $n\in\NN$}.$$
$\dim(E)$ is called the index of the dominated splitting.
 
 The following result is well known, which tells us that the dominated bundles are invariant under the diffeomorphisms in the centralizer.  
 
 \begin{Proposition}[Lemma 13 in \cite{DWX}]~\label{p.preserve-domination}
 	Let $f\in\diff^1(M)$ admit a dominated splitting of the form $TM=E\oplus F$.  Then for any $g\in\cZ^1(f)$, one has $Dg(E)=E$ and $Dg(F)=F$.
 \end{Proposition}

\subsection{Dynamical coherence}
For a partially hyperbolic diffeomorphism $f$, by \cite{HPS}, there always exist $f$-invariant foliations $\cF^s$ and $\cF^u$ tangent to the bundles $E^{s}$ and $E^u$ respectively and such foliations are unique.  
 $f$ is called \emph{dynamically coherent} if there exist $f$-invariant foliations $\cF^{cs}$ and $\cF^{cu}$ tangent to  $E^{cs}:=E^{s}\oplus E^c$ and $E^{cu}:=E^c\oplus E^u$ respectively.  By taking the intersection of $\cF^{cs}$ and $\cF^{cu}$, one gets an invariant foliation $\cF^c$ tangent to $E^c.$
The dynamical coherence of a partially hyperbolic DA-diffeomorphism  on $\TT^3$ has been substantially investigated, see for instance~\cite{BBI1, BI, Po, HaPo0, FPS}.

Two transverse foliations $\cF,\cG$ on $\RR^3$ have \emph{global product structure} if for any $x,y\in\RR^3$, the leaf $\cF(x)$ intersects the leaf $\cG(y)$ into a unique point. Foliation $\cF$ on $\RR^3$ is \emph{quasi-isometric} if there exist $a,b>0$ such that for any $x\in\RR^3$ and $y\in\cF(x)$, one has $\ud_\cF(x,y)\leq a\cdot \ud(x,y)+b$, where $\ud_\cF(\cdot,\cdot)$ denotes the distance on the leaves of $\cF$ and $\ud(\cdot,\cdot)$ denotes the Euclidean distance. 

The following result gives the dynamical coherence of a partially hyperbolic hyperbolic DA-diffeomorphism and further geometrical properties of the invariant foliations.
\begin{Theorem} [\cite{BI,Ha,Po}]~\label{thm.dynamical-coherence-unique}
 Let $f\in\diff^1(\TT^3)$ be a partially hyperbolic diffeomorphism with the partially hyperbolic splitting $T\TT^3=E^{s}\oplus E^c\oplus E^{u}$. Assume that $L_f$ is Anosov, 
	then one has the following 
	\begin{itemize}
		\item $f$ has   unique foliations $\cF^{cs}$ and  $\cF^{cu}$ tangent to $E^{s}\oplus E^c$ and  $E^c\oplus E^{u}$ respectively;
		\item  the lifts of the foliations $\cF^{cs}$ and $\cF^u$ to $\RR^3$ have global product structure;
		\item the lifts of the foliations $\cF^s,\cF^c,\cF^u$ to $\RR^3$ are quasi-isometric;
		\item each leaf of $\cF^c$ is dense in $\TT^3$;
		\item $L_f$ has simple spectra.
	\end{itemize}
\end{Theorem}
\begin{Remark}
Key Lemma 2.1	in \cite{BI} gives the existence of 2-dimensional  foliations transverse to $E^{u}$ and $E^{s}$ respectively which is exactly the assumption of  Theorem A in \cite{Po}.
\end{Remark}

\begin{Notation}
	Throughout this paper, for any foliation $\cF$ on $\TT^3$, we will denote by $\tilde{\cF}$ 
	the lift of $\cF$ to $\RR^3.$ We denote $\ud_{\tilde{\cF}}(\cdot,\cdot)$ the distance in a $\tilde{\cF}$-leaf, and 
	$$
	\tilde{\cF}_r(x)=
	\big\{y\in\tilde{\cF}(x): \ud_{\tilde{\cF}}(x,y)<r\big\}.
	$$ 
	And we assume the center Lyapunov exponent of $L_f$ is larger than zero, i.e. the stable dimension of $L_f$ as an Anosov diffeomorphism is $1$. Otherwise, we only need to consider $f^{-1}$. 
\end{Notation}

As a consequence of the global product structure for the lifted foliations, one has the following result whose proof can be found in~\cite[Proposition 6.8]{Po}. 

\begin{Corollary}~\label{co:quantitative-product-structure}
	Let $f$ be as in the assumption of Theorem~\ref{thm.dynamical-coherence-unique}. Then there exists a constant $K>0$ such that for any $r>0$, any $x\in \RR^3$, any  $y\in\tilde{\cF}_{r}^{u}(x)$ and any $w\in\tilde{\cF}^{cs}(x)$, one has that  $\tilde{\cF}^u_{r+K}(w)\cap\tilde{\cF}^{cs}(y)\neq\emptyset$.
	
	The analogous result with respect to the strong stable and center unstable foliations  holds.
\end{Corollary}
 Combining with the result from the previous section, one has the following corollary:
 \begin{Corollary}~\label{co:preserve-foliations}
 Let 	$f$ be a $C^1$-partially hyperbolic diffeomorphism homotopic to a linear Anosov.
 Assume that $f$ has the splitting of the form $TM=E^s\oplus E^c\oplus E^u,$
  then for any $g\in\cZ^1(f)$, each invariant foliation  $\cF^*$ of $f$ is invariant under $g$ for $*=s,cs,c,cu,u.$
 	\end{Corollary}
 \proof
 By the classical stable manifold theorem and Theorem~\ref{thm.dynamical-coherence-unique}, there exist unique invariant foliations $\cF^s$, $\cF^u$, $\cF^{cs}$ and $\cF^{cu}$ tangent to $E^s$, $E^u$, $E^{cs}$, and $E^{cu}$
 respectively. For any $g\in\cZ^1(f)$, by Proposition~\ref{p.preserve-domination}, one has that $g(\cF^*)$ is an $f$-invariant foliation tangent to $E^*$ for $*=s,u,cs,cu$; therefore, one has $g(\cF^*)=\cF^*$ for $*=s,u,cs,cu$. 
 Finally, we have
 $$
 g(\cF^c)=g(\cF^{cs})\cap g(\cF^{cs})=\cF^{cs}\cap\cF^{cs}=\cF^c.
 $$
 \endproof
 
 For a diffeomorphism  on the torus, if its linear part is Anosov, then it is semi-conjugate to its linear part. 
 \begin{Theorem}[\cite{F, Wa}]\label{thm.unique-semi-conjugacy}
Let $f\in\diff^1(\TT^d)$ and assume that $L_f$ is  Anosov. Consider a lift   $F$   of $f$ to the universal cover $\RR^d$,  then there exists a unique continuous surjective map $H: \RR^d\to\RR^d$ such that 
 \begin{itemize}
 	\item $H\circ F=L_f\circ H$;
 	\item  $H(x+z)=H(x)+z$, for any $z\in\ZZ^d$ and any $x\in\RR^d$.
 \end{itemize}
\end{Theorem}
As a consequence, one has the following:
\begin{Corollary}~\label{co:linearizable-simultaneously}
	Let $f\in\diff^1(\TT^d)$ whose linear part $L_f$ is Anosov and $F$ be a lift of $f$ to $\RR^d$. Then there exists a continuous surjective map $H:\RR^d\to\RR^d$ such that 
	\begin{itemize}
		\item $H\circ F=L_f\circ H$;
		\item $H(x+z)=H(x)+z$ for any $x\in\RR^d$ and any $z\in\ZZ^d$;
		\item for  any $g\in\cZ^1(f)$ and any lift   $G$ of $g$ to $\RR^d$, if $F\circ G=G\circ F,$ then $H\circ G=L_g\circ H.$
	\end{itemize} 
\end{Corollary}
\proof
Let $H:\RR^d\to\RR^d$ be the continuous surjective map given by Theorem~\ref{thm.unique-semi-conjugacy} such that $H\circ F=L_f\circ H$ and $H-\Id_{\RR^d}$ is $\ZZ^d$-periodic. 
 Consider the map $\wh{H}=L_g^{-1}\circ H\circ G$ which satisfies that $\wh H-\Id_{\RR^d}$ is  $\ZZ^d$-periodic. Then one has 
$$\wh H\circ F=L_g^{-1}\circ H\circ G\circ F=L_g^{-1}\circ H\circ F\circ G=L_g^{-1}\circ L_f\circ H\circ G=L_f\circ L_g^{-1} \circ H\circ G=L_f\circ \wh H.$$
By the uniqueness property in Theorem~\ref{thm.unique-semi-conjugacy}, one has $\wh H=H$ which gives $H\circ G=L_g\circ H$
\endproof 
Furthermore, the semi-conjugation preserves certain foliations.
\begin{Theorem}[\cite{HaPo0,Ha, Po,U}]~\label{thm.semi-conjugacy}
	Let $f$ be a $C^1$-partially hyperbolic diffeomorphism on $\TT^3$ which is homotopic to an Anosov automorphism $L_f$ with two positive Lyapunov exponents. Denote  $\cF^*$ and $\mathcal{W}^*$ the foliations of $f$ and $L_f$ respectively, for $*=s,u,cu,cs,c$. Let $h$ be the semi-conjugacy between $f$ and $L_f$, in formula: $L_f\circ h=h\circ  f.$
	Then one has the following properties:
	\begin{enumerate}
		\item For any $x\in\TT^3$ and $*=cu,cs,c$, one has $h(\cF^*(x))=\cW^*(h(x)).$ 
		\item  For any $x\in\TT^3,$ the map $h|_{\cF^{s}(x)}:\cF^{s}(x)\rightarrow \cW^{s}(h(x))$ is a homeomorphism. 
		\item  For any $x\in\TT^3$,  the pre-image $h^{-1}(h(x))$ is a segment (could be trivial) contained in $\cF^c(x)$. In particular, for any center leaf $\cW^c(y)$, there exists at most countable many points whose pre-image under $h$ are non-trivial center segments. 
	\end{enumerate}
\end{Theorem}

\begin{Remark}~\label{r.ures}
	\begin{enumerate}
		\item  This result is obtained in \cite{Ha, U} assuming  absolute partial hyperbolicity and it is extended to general partially hyperbolic setting in \cite[Appendix A]{Po} (see also \cite[Section 3]{HaPo0}).
		\item The last item implies that $f$ must have fixed points.
		\item Each center leaf of $f$ is dense in $\TT^3$. 
	\end{enumerate}
\end{Remark}

\subsection{Centralizer of   linear Anosov automorphisms}
The following result comes from \cite{AP} (see also \cite{W1}).
\begin{Theorem}[\cite{AP}]~\label{thm.adler-palais}
	Let $L$ be a linear Anosov map on $\TT^3$ and $h$ be a homeomorphism on $\TT^3$. 
	If $h\circ L=L\circ h$, then $h$ is affine. 
\end{Theorem}


In particular, Adler-Palais' result implies that for each Anosov diffeomorphism $f$ on $\TT^3$, there exists a homeomorphism $h$ on $\TT^3$ such that for any $g\in\cZ^0(f)$, one has that $h\circ g\circ h^{-1}$ is affine. Corollary~\ref{co:linearizable-simultaneously} tells us that such results hold  for semi-conjugacy to Anosov case.

The following result gives the rank of the group of linear automorphisms commuting with an Anosov automorphism,  and it comes from \cite{KKS,DWX}.
\begin{Lemma}[Proposition 3.7 in \cite{KKS} and Lemma 16 in \cite{DWX}]\label{l.rank-of-centralizer}
	Consider a matrix $L\in SL(n,\ZZ)$ whose characteristic polynomial is irreducible over $\ZZ$. 
	Then the group $\cG(L)=\big\{L_1\in SL(n,\ZZ): LL_1=L_1L\big\}$ is abelian. Moreover, $\cG(L)$ is virtually $\ZZ^{r+c-1}$, where $r$ is the number of real eigenvalues of $L$ and $2c$ is the number of the complex eigenvalues of $L$.
\end{Lemma}
In our paper, all the linear Anosov maps we consider have real simple spectra. It is easy to see that these linear Anosov maps are irreducible in the sense that their characteristic polynomials are irreducible over $\ZZ$ \footnote{Otherwise,  one should have $\pm 1$ as the eigenvalue, contradicting to the hyperbolicity of the maps.}.Therefore, their centralizers are  virtually $\ZZ^2$ by Lemma~\ref{l.rank-of-centralizer}.  Moreover,  the non-trivial elements in the centralizer is also Anosov.
\begin{Lemma}~\label{l.centralizer-of-linear-Anosov}
	Let $A\in SL(3,\ZZ)$ be an Anosov automorphism. For any $B\in SL(3,\ZZ)$, if $AB=BA$,  then either $B=\Id$ or $B$ is Anosov.
\end{Lemma}
\proof
If $B\in SL(3,\ZZ)$ has eigenvalues of modulus $1$,  then $1$ is an eigenvalue of $B$ since $\det(B)=1$. And  there exists a rational vector  $0\neq v$ such that $Bv=v$. For any eigenvector $w$ of $B$ with respective to $1$, one has    that $BAw=ABw=Aw$, which implies that the eigenvector space of $B$ with respect to $1$ is invariant under $A$; in particular  the rational vector $Av$ is  also an eigenvector  of $B$ with respect to the eigenvalue $1$ and is not collinear to $v$,  since $A$ is Anosov and $v$ is rational.   Then the 2-dimensional linear space generated by $Av$ and $v$ 
has rational slope and 	is   contained in the eigenspace of $B$ with respect to the eigenvalue $1$.  Once again, as $A$ is Anosov (in particular irreducible), the linear  space generated by $Av$ and $v$ is not $A$-invariant which implies that each vector in $\RR^3$ is in the eigenspace of $B$ with respect to eigenvalue $1$. Hence $B$ is identity.
\endproof

\subsection{Regularity}

Now, we collect some regularity lemmas showing that if a homeomorphism is differentiable along pairs of transverse foliations up to certain order, then the homeomorphism is differentiable. 
\begin{Lemma}[\cite{J}]~\label{l.journe}
	Let $M$ be a closed manifold and $h$ be a homeomorphism on $M$. Assume that there exist two  transverse continuous  foliations $\cF$ and $\cG$ on $M$ with $C^r$-leaves, and 
	$h$ is uniformly $C^r$ when restricted to leaves of $\cF$ and $\cG$, then $h$ is $C^{r-\e}$ for any $\e>0$.
\end{Lemma}

\begin{Remark}
	If $r$ is not an integer, then $h$ is $C^r$. If $r$ is an integer, then $h$ is $C^{r-1+{\rm Lip}}$. 
\end{Remark}

One says that a foliation $\cF$ with $C^1$-leaves is  \emph{expanding} for $f\in\diff^1(M)$ if $\cF$ is $f$-invariant and there exists $N\in\NN$ such that $\|Df^{-N}|_{T_x\cF(x)}\|\leq\frac{1}{2},$  for any $x\in M$.
\begin{Lemma}[Lemma 2.4 in \cite{G}]~\label{l.bootstrap}
	Let $f,g$ be two $C^r$-diffeomorphisms on a closed manifold $M.$  Let $\cF, \cG$ be one dimensional expanding foliations with $C^r$-leaves for $f$ and $g$ respectively. Assume that there exists a homeomorphism $h$ on $M$ such that 
	\begin{itemize}
		\item  $h\circ f=g\circ h$ and  $h(\cF)=\cG$;
		\item  $h$ and its inverse are uniformly $C^1$ along the leaves  of $\cF$ and $\cG$ respectively;
	\end{itemize} 
	then $h$ is uniformly $C^r$ along the leaves of $\cF$ and $h^{-1}$ is uniformly $C^r$ along the leaves of $\cG.$ 
\end{Lemma}

\subsection{Accessibility}
Given a partially hyperbolic diffeomorphism $f\in\diff^1(M)$  and a point $x\in M$, the \emph{accessible class $\Acc(x)$ of $x$} is defined as the set of points which can be joined to $x$ by paths which are concatenations of paths in a strong stable or strong unstable manifold. By definition, $\Acc(x)$ is saturated by strong stable and strong unstable leaves. One says that $f$ is \emph{accessible} if  any two points $x,y\in M$ can be connected by   a path $\gamma$ which is a concatenation of paths in strong stable or stable unstable manifolds of $f$, in other words, the accessible class of a point is the whole manifold. 

For a partially hyperbolic diffeomorphism $f$, the bundles $E^s$ and $E^u$ are \emph{jointly integrable}, if there exists an $f$-invariant foliation tangent to $E^s\oplus E^u$ everywhere. In this case, we    call $f$ is su-integrable or su-jointly integrable.

It has been proved in \cite{Di} that if $f$ is accessible, then $E^s\oplus E^u$ is not jointly integrable. Moreover, if $f$ has one-dimensional center, then there exists a point $x\in M$ has the local accessibility property as following.

\begin{Lemma}[\cite{Di,HHU}]~\label{l.twist-in-center}
	Let  $f$ be a $C^1$-partially hyperbolic diffeomorphism on $M$. If $f$ is accessible and the center bundle is one dimensional, then there exist $r_0>0, r_1>0$ which can be arbitrarily small, and $x\in M$ such that  any center curve $I_{r_1}^c(x)$ centered at $x$ of radius $r_1$, there exist $x^s,x^{su}, x^{sus}\in M$ and $x^c\in I_{r_1}^c(x)$ such that 
	\begin{itemize}
		\item $x^s\in \cF^{s}_{r_0}(x)$ and  $x^{su}\in\cF^{u}_{r_0}(x^s)$ 
		\item $x^{sus}\in\cF^s_{r_0}(x^{su})$  and $x^c\in\cF^u_{r_0}(x^{sus})$,
	\end{itemize}
	where $\cF^*_{r}(z)$ denotes the $r$-neighborhood of $z$ in the leaf $\cF^*(z)$ for $*=s,u$. 
	
Moreover,   let  $I^c(x,x^c)$ denote the set   of all points located between $x$ and $x^c$ in $I_{r_1}^c(x)$, then  $\Acc(x)$ contains an open set $U$ close to $x$, i.e.
	$$
	U~\subseteq~ \bigcup_{y\in I^c(x,x^c)}
	\bigcup_{z\in\cF^s_{\it loc}(y)}\cF^u_{\it loc}(z) 
    ~\subseteq~ \Acc(x).
	$$
	Each point in  $U$ can be connected to $x$ by a local $su$-path contained in a small neighborhood of $x$.
\end{Lemma}

If a diffeomorphism $f:M\rightarrow M$ is partially hyperbolic, and $\pi:\tilde M\rightarrow M$ is a covering map, then any lift $\tilde f:\tilde M\rightarrow \tilde M$ is also partially hyperbolic. And the partially hyperbolic splitting on $\tilde M$ is defined by pulling back the splitting on $M$:
$$
T\tilde M=\pi^*(E^s)\oplus \pi^*(E^c)\oplus \pi^*(E^u).
$$
The following result tells us that the accessibility is preserved under lifts of the manifold.

\begin{Lemma}\label{l.lift-preserve-accessibility}
	Let  $f\in\diff^1(M)$ be an accessible partially hyperbolic diffeomorphism, and assume that  the center bundle is one dimensional. Consider a covering map  $\pi: \tilde M\to M$ from a connected manifold $\tilde M$ to $M$, and a lift $\tilde f$ of $f$ to  $\tilde M$, then $\tilde f$ is accessible.
\end{Lemma}
\proof
Notice that the lift of the strong stable and unstable foliations to $\tilde M$ are the strong stable and unstable foliations of $\tilde f$. Since $f$ is accessible, Lemma~\ref{l.twist-in-center} shows that there exists $x\in M$ which has the local accessibility property.  

This implies for every $\tilde x\in\pi^{-1}(x)$, the accessibility class $\Acc(\tilde{x})$ with respect to $\tilde f$ contains an open set close to $\tilde x$. If an accessible class contains an open set, then it is open. 
Thus $\Acc(\tilde x)$ is open for every $\tilde x\in\pi^{-1}(x)$. 

On the other hand, for every $\tilde y\in\tilde M$, since $\pi(y)\in M$ can be connected to $x$ by an $su$-path, $\tilde y$ can be connected to a point $\tilde x\in\pi^{-1}(x)$ by an $su$-path. This implies
$$
\tilde M=\cup_{\tilde x\in\pi^{-1}(x)}\Acc(\tilde{x}).
$$
Since $\tilde M$ is connected and each $\Acc(\tilde{x})$ is open, we must have $\tilde M=\Acc(\tilde{x})$ for every $\tilde x\in\pi^{-1}(x)$. Thus $\tilde f$ is accessible.
\endproof


The following result gives equivalence conditions for the joint integrability of strong stable and unstable distributions. 
\begin{Theorem}[Theorem 1.1 in \cite{GS} and \cite{HaS}]\label{thm.DA-is-ergodic}
	Let $f$ be a $C^r$ ($r>1$) partially hyperbolic diffeomorphism on $\TT^3$ whose linear part $L_f$ is Anosov.  The followings are equivalent:
	\begin{itemize}
		\item strong stable and unstable distributions of $f$ are jointly integrable;
		\item  each periodic orbit of $f$ has same center Lyapunov exponent as $L_f$, and $f$ is Anosov;
		\item  $f$ is not accessible.
	\end{itemize}
\end{Theorem}
\begin{Remark}
	\begin{itemize}
		\item  Under volume preserving assumption, A. Hammerlindl and R. Ures \cite{HaU} proved that  the first and   third items are equivalent; in particular,  $f$ is topological Anosov;
		\item The equivalence of the first and second items is obtained in~\cite{GS} under  volume preserving setting. Then the volume preserving condition is removed by \cite{HaS} and the third equivalent item is obtained in ~\cite{HaS}.
	\end{itemize}
	
\end{Remark}

\subsection{Equivalent conditions for su-integrability of an Anosov map }

In this part, we collect the consequences of su joint integrability for Anosov diffeomorphisms on $\TT^3$, which is proved in~\cite{GS} and \cite{GRZ}.

\begin{Theorem}[Theorem 5.1 in \cite{GS}]~\label{thm.gan-shi}
	Let $f$ be a $C^r$ ($r>1$) partially hyperbolic and Anosov diffeomorphism on $\TT^3$. Let $h$ be the conjugacy between $f$ and $L_f$. 
	Then the followings are equivalent:
	\begin{itemize}
		\item $f$ is $su$-integrable;
		\item  $f$ is not accessible;
		\item $h$ preserves the strong stable and strong unstable foliations;
		\item   the center Lyapunov exponent of any periodic point $p$ of $f$ coincides with the center Lyapunov exponent of $L_f$;
		\item  $h$ is differentiable along the center leaves of $f$.
	\end{itemize}
\end{Theorem}
\begin{Remark}\label{r.gan-shi}
	When the conjugacy preserves the strong foliations, one can show that $h$ and $h^{-1}$ is uniformly H\"older continuous along the   leaves of strong foliations (see for instance Lemma 2.3 in~\cite{GS}).
\end{Remark}

Now, we state the following result which is essential Proposition 4.1  in \cite{GS}. For completeness, we will give the proof in Appendix~\ref{s.appendix}.
\begin{Theorem}~\label{thm.generalized-gan-shi-on-equi-center-exponent}
	Let $g$ be a  $C^r$ ($r>1$)  Anosov diffeomorphism on $\TT^3$ and let $h\in\homeo(\TT^3)$ such that $h\circ g=L_g\circ h$. Assume that 
	\begin{itemize}
		\item there exists $Dg$-invariant continuous splitting $E^s\oplus E^c\oplus E^u$;
		\item  $g$ is uniformly contracting along $E^s$ and is uniformly expanding along $E^c\oplus E^u$;
		\item   there exist $g$-invariant foliations $\cF^c$, $\cF^u$ and $\cF^{su}$ tangent to $E^c$,  $E^u$ and $E^s\oplus E^u$ respectively;
		\item  $L_g$ is partially hyperbolic;
		\item  $h$ sends  $\cF^c, \cF^u$ to the center, strong unstable  foliations of $L_g$ respectively;
		\item the holonomy map given by $\cF^u$ restricted to each unstable leaf between two local plaques tangent to $E^c$ at a uniform bounded distance is uniformly $C^1$;
	\end{itemize}
	Then the Lyapunov exponent along $E^c$ of  any periodic point $p$ is same as the center Lyapunov exponent of $L_g$, and $h$ is uniformly $C^1$ along the leaves of $\cF^c$.
\end{Theorem}
\begin{Remark}
	In the statement of Theorem~\ref{thm.generalized-gan-shi-on-equi-center-exponent}, we do not assume the splitting $E^c\oplus E^u$ is dominated. 
\end{Remark}

\section{Centralizer of partially hyperbolic DA-diffeomorphism:\\ Proof of Theorems~\ref{thm.main-thm} and \ref{thm.main-thm-low-regularity}}~\label{s.proof-of-thms}
In this section, we give the proof of our main theorems. 
For $f\in\diff^r(\TT^3)$, the linearization   of the centralizer of $f$ is the group $\{L_g: g\in\cZ^r(f) \}$. The proof is proceeded according to the linear part of the centralizer of    partially hyperbolic DA-diffeomorphisms. We will first discuss the case where the linearization  of the centralizer is virtually  $<L_f>$. Then we discuss the case where the group $\{L_g: g\in\cZ^r(f) \}/<L_f>$ is not virtually trivial.  

\subsection{Preliminary lemmas}
The lifts of two commutable diffeomorphisms may  not be  commutable. The following result tells us that the lifts of the centralizer of partially hyperbolic DA-diffeomorphisms, up to finite iterates,  are still in the centralizer of the lifted diffeomorphism.

\begin{Lemma}~\label{l.commutable-up-toiterations}
	Let $f$ be a $C^1$-partially hyperbolic diffeomorphism on $\TT^3$ whose linear part $L_f$ is Anosov, and let  $F$ be a lift of $f$ to $\RR^3$. 
	
	Then for any $g\in\cZ^1(f)$, there exist an integer $0<l\leq |\det(L_f-\Id_{\RR^3})|$ and a lift $\wh G$ of $g^l$ such that $F\circ \wh G=\wh G\circ F.$  Furthermore, if $g$ is homotopic to identity,  then  $l$ can be chosen as a factor of $|\det(L_f-\Id_{\RR^3})|$.
\end{Lemma}
\proof
Let $p_1,\cdots,p_k\in\TT^3$ be all the fixed points of $L_f.$ It is classical that $k=|\det(L_f-\Id_{\RR^3})|$.     Let $H:\RR^3\to\RR^3$ be the semi-conjuacy between $F$ and $L_f$ given by Theorem~\ref{thm.unique-semi-conjugacy}, and let  $h:\TT^3\to\TT^3$ be the map induced by $H$. 
 Then the set of fixed points of $f$ is contained in $\cup_{i=1}^k h^{-1}(p_i)$ and each $h^{-1}(p_i)$ is $f$-invariant. By Theorem~\ref{thm.semi-conjugacy}, each $h^{-1}(p_i)$ is a compact center segment (could be trivial). By Brouwer fixed points theorem, $f$ has fixed points in each $h^{-1}(p_i)$.  Let $I_i\subset h^{-1}(p_i)$ be the shortest connected and compact center segment (could be trivial) containing all fixed points of $f$ in $h^{-1}(p_i)$. Then the two endpoints of $I_i$ are fixed points of $f$. 

Let  $\pi:\RR^3\to\TT^3$ be the canonical covering map. Without loss of generality, one can assume that $p_1$ is the projection of $0\in\RR^3$ under $\pi$, i.e. $\pi(0)=p_1.$

Since $g\in\cZ^1(f)$, by Corollary~\ref{co:preserve-foliations}, $g$ preserves the center foliation of $f$. As the set of fixed points of $f$ is $g$-invariant, for each $i\in\{1,\cdots,k \}$, there exists $j\in\{1,\cdots, k \}$ such that $g(I_i)=I_j$ which defines a permutation on $\{1,\cdots, k\}.$  Therefore, there exists $0<l\leq k$ such that $g^{l}(I_1)=I_1$. If $g^l$ preserves the orientation of the center bundle,  then the two endpoints of $I_1$ are the fixed points of $f$ and $g^l$. If $g^l$ reverses the orientation of the center bundle, then $g^l$ has a unique fixed point in $I_1$ which is also a fixed point of $f$ since $g^l(I_1)=f(I_1)=I_1$ and $g\in\cZ^1(f)$. To summarize, $f$ and $g^l$ have a common fixed point $q_1\in I_1$. Notice that  $\wh q_1=H^{-1}(0)\cap \pi^{-1}(q_1)$ is a fixed point  of $F$.  Since $q_1$ is a fixed point of $g^l$, there exists a lift $\wh G$ of $g^l$ such that $\wh G(\wh q_1)=\wh q_1.$   Observe that    $F\circ \wh G\circ F^{-1}\circ \wh G^{-1}$ is a lift of identity map on $\TT^3$ and has a fixed point $\wh q_1$, hence $F\circ \wh G=\wh G\circ F.$ 
\medskip

Now, we assume that $g$ is homotopic to identity. Let $G$ be a lift of $g$ to $\RR^3.$ Since $f\circ g=g\circ f$,  there exists  $n\in\ZZ^3$ such that $F\circ G=G\circ F+n.$ Since $L_f$ is Anosov, the linear map  $L_f-\Id_{\RR^3}$ is  invertible. Let $m=(L_f-\Id_{\RR^3})^{-1}n\in\QQ^3$, then there exists an integer $l>0$ which is a factor of $|\det(L_f-\Id_{\RR^3})|$ such that  $l\cdot m\in\ZZ^3$. Since $g$ is homotopic to identity, then  $F\circ G^l=G^l\circ F+ln$. Now, let $\wh G=G^l-l\cdot m$ which is a lift of $g^l$,  and one has the following 
\begin{align*}
F\circ\wh G=F\circ(G^l-lm)=F\circ G^l-L_f(lm)=F\circ G^l-ln-lm=G^l\circ F-lm=\wh G\circ F,
\end{align*}  which ends the proof.
\endproof

The following result discusses the existence of  common fixed points for lifted dynamics.
\begin{Lemma}~\label{l.common-fixed-point-for-lifts}
		Let $f$ be a $C^1$-partially hyperbolic diffeomorphism on $\TT^3$ whose linear part $L_f$ is Anosov, and let  $g\in\cZ^1(f)$. Assume that there exist a lift $F$ of $f$ to $\RR^3$ and a lift $G$ of $g$ to  $\RR^3$ such that $F\circ G=G\circ F$, then $F$ and $G$ have a common fixed point, that is,  there exists $p\in\RR^3$ such that $F(p)=G(p)=p$.   
\end{Lemma}
\proof
By Corollary~\ref{co:preserve-foliations}, the center foliation of $F$ is  $G$-invariant.
Let $H:\RR^3\to\RR^3$ be the continuous surjective map given by Corollary~\ref{co:linearizable-simultaneously} such that 
\begin{itemize}
	\item $H\circ F= L_f\circ H$ and $H\circ G=L_g\circ H;$
	\item  $H-\Id_{\RR^3}$ is $\ZZ^3$-periodic.
\end{itemize}
As $L_f$ and $L_g$ have a unique fixed point $0\in\RR^3$, then all the fixed points of $F$ and $G$ are contained in $H^{-1}(0)$.  By Theorem~\ref{thm.semi-conjugacy}, $H^{-1}(0)$ is a compact  and $F$-invariant center segment. 

As $F$ commutes with $G$, the set of fixed points of $F$ is $G$-invariant and vice versa.
If $F$ reserves the orientation of the center bundle, by the fact that $\dim(E^c)=1$,   $F$ has a unique fixed point in $H^{-1}(0)$ which is also  a fixed point of $G$. If $G$ reverses the orientation of the center bundle, one concludes analogously. 
 If  $F$ and $G$   preserve  the orientation of the center bundle, the  endpoints of $H^{-1}(0)$ are the fixed points of $F$ and $G$, proving the existence of common fixed points.
\endproof 
\subsection{The linearization of the centralizer is virtually $\ZZ$: virtually trivial case}
In this part, we discuss the case where the group $\{L_g: g\in\cZ^r(f) \}$ is virtually $\ZZ,$ and we show in this case that  the centralizer of $f$ is virtually trivial.  

The main purpose of this section is to prove the following result. 

\begin{Theorem}~\label{thm.virtually-trivial-case-homotopic-to-id}
	Let $f$ be a $C^r$ ($r>1$) partially hyperbolic diffeomorphism on $\TT^3$ whose linear part $L_f$ is Anosov. 
	Then one has that 
	\begin{itemize}
		\item for each $g\in\cZ^r(f)$ which is homotopic to identity, there exists an integer $l=l_g$ which is a factor of $|\det(L_f-\Id_{\RR^3})|$ such that $g^l=\Id_{\TT^3}.$
		\item 	$$\#\big\{g\in\cZ^r(f): \textrm{$g$ is homotopic to $\Id_{\TT^3}$} \big\}\leq |\det(L_f-\Id_{\RR^3})|.$$
	\end{itemize}

\end{Theorem} 

Before giving the proof of Theorem~\ref{thm.virtually-trivial-case-homotopic-to-id}, we need to make some preparations.
We first show that if $g\in\cZ^r(f)$ is homotopic to identity and a lift of $g$ admits a  fixed point on  $\RR^3$, then $g$ is identity.

\begin{Proposition}~\label{p.center-fixing}
	Let $f\in\diff^r(\TT^3)$ ($r>1$) be a partially hyperbolic diffeomorphism and $g\in\cZ^r(f)$. 
	Assume that 
	\begin{itemize}
		\item the linear part $L_f$ of $f$ is Anosov;
		\item  there exists a  lift  $G:\RR^3\to\RR^3$ of  $g$   with the following properties:
	\begin{itemize}
		\item[--] $G(x+n)=G(x)+n$ for any $x\in\RR^3$ and $n\in\ZZ^3$;
		\item[--] $G$ admits a fixed point   $q\in\RR^3$;
	\end{itemize} 
	\end{itemize}
then $G=\Id_{\RR^3}.$
	\end{Proposition}
\proof
By Corollary~\ref{co:preserve-foliations} and the fact that $G-\Id_{\RR^3}$ is $\ZZ^3$-periodic, one has $$G\big(\tilde\cF^c(q+n)\big)=\tilde\cF^c(q+n)=\tilde\cF^c(q)+n \textrm{ for any $n\in\ZZ^3.$} $$
By the third  item of Remark~\ref{r.ures}, the set $\big\{\tilde\cF^c(q+n)\big\}_{n\in\ZZ^3}$ is dense in $\RR^3$.

The following claim tells us that $G$ is center fixing.
\begin{Claim-numbered}~\label{cl.center-fixing}
	For any $x\in\RR^3$, the center leaf $\tilde\cF^c(x)$ is fixed by $G$.
\end{Claim-numbered}
\proof[Proof of the claim]
For any $x\in\RR^3$, there exists a sequence of points $x_k\in\big\{\tilde\cF^c(q+n)\big\}_{n\in\ZZ^3}$ such that $x_k$ converges to $x$.
Since $G-\Id_{\RR^3}$ is $\ZZ^3$-periodic, there exists $\ell_0>0$ such that $\ud(G(y),y)\leq\ell_0$  for any $y\in\RR^3$.  Since the center foliation $\tilde{\cF}^c$ is quasi-isometric, there exist $a,b>0$ such that for any $x,y\in\RR^3$ with $x\in\tilde{\cF}^c(y)$, one has $\ud_{\tilde\cF^c}(x,y)\leq a\cdot \ud(x,y)+b.$
Since the center leaf $\tilde\cF^c(x_k)$ is $G$-invariant, one has
$$\ud_{\tilde\cF^c}(G(x_k),x_k)\leq a\ud(G(x_k),x_k)+b\leq a\ell_0+b.$$
Let $\ell_1=a\ell_0+b.$ By the continuity of the center foliation, $\tilde{\cF}_{\ell_1}^c(x_k)$ converges to $\tilde{\cF}^c_{\ell_1}(x)$. By the continuity of $G$ and the fact that $G(x_k)\in\tilde{\cF}_{\ell_1}^c(x_k)$, one has   $G(x)\in\tilde{\cF}^c_{\ell_1}(x)$ proving that the center leaf $\tilde{\cF}^c(x)$ is fixed by $G$.
\endproof 
\begin{Claim-numbered}\label{Cl: Amie}
	For any fixed point $x_0$ of $G$, one has  
$$G|_{\tilde\cF^s(x_0)\cup\tilde\cF^u(x_0)}=\Id|_{\tilde\cF^s(x_0)\cup\tilde\cF^u(x_0)}.$$
\end{Claim-numbered}

\proof[Proof of the claim]
Let $x_0$ be a fixed point of $G$. By Corollary~\ref{co:preserve-foliations}, one has $G(\tilde{\cF}^s(x_0))=\tilde{\cF}^s(x_0)$ and $G(\tilde\cF^u(x_0))=\tilde\cF^u(x_0)$.  One only needs to show that $G$ restricted to $\tilde\cF^s(x_0)$ is identity, and the case for the strong unstable manifold is analogous. Now, consider the center stable leaf $\tilde\cF^{cs}(x_0)$. By Theorem~\ref{thm.dynamical-coherence-unique}, the foliations $\tilde\cF^{cu}$ and $\tilde{\cF}^{s}$ have global product structure,  
which implies that for any  point $y\in\tilde\cF^{cs}(x_0)$,  the leaf $\tilde\cF^c(y)$ intersects $\tilde\cF^s(x_0)$ into a unique point. For any $y\in\tilde\cF^s(x_0)$, by Claim~\ref{cl.center-fixing} and the fact that $G(\tilde\cF^s(x_0))=\tilde\cF^s(x_0),$  one has $\{G(y)\}=G\big(\tilde\cF^c(y)\cap\tilde\cF^s(x_0)\big)=\tilde\cF^c(y)\cap\tilde\cF^s(x_0)=\{y\}.$ 
\endproof

Now, we show that $G$ is identity.
As $\Acc(q)$   is saturated by strong stable and strong unstable leaves and $q$ is a fixed point of $G$,  by Claim~\ref{Cl: Amie},   
the map $G$ coincides with identity on $\Acc(q)$.
There are two cases to discuss according to  accessible property.

If $f$ is accessible, by Lemma~\ref{l.lift-preserve-accessibility},  each lift of $f$ to the universal cover is also accessible, hence $\Acc(q)=\RR^3$ which implies $G=\Id_{\RR^3}.$

If $f$ is not accessible, by Theorem~\ref{thm.DA-is-ergodic}, $f$ is   Anosov. 
Consider a projection $p$ of the fixed point $G$ on $\TT^3$, then $g$ coincides with identity on the union of the strong stable and unstable manifolds of $p$.  As $f$ is Anosov, then the union of the strong stable and unstable manifolds of $p$ is dense in $\TT^3$, hence $g=\Id_{\TT^3}$ which in return implies $G=\Id_{\RR^3}$ since $G$ has fixed points. Now the proof of Proposition~\ref{p.center-fixing} is completed.
\endproof

As a consequence, one has the following corollary. 
\begin{Corollary}~\label{thm.homotopic-to-id}
	Let $f\in\diff^r(\TT^3)$ ($r>1$) be a  partially hyperbolic diffeomorphism  whose linear part $L_f$ is Anosov,  and let $g\in\cZ^r(f)$.
	If $g$ is homotopic to identity, then there exists  $l\in\NN$ which is a factor of $|\det(L_f-\Id_{\RR^3})|$ such that $g^l=\Id_{\TT^3}.$ 
\end{Corollary} 
\proof 
Consider a lift of $F$ of $f$ to $\RR^3$. Let $g\in\cZ^r(f)$ be a diffeomorphism homotopic to identity. By Lemma~\ref{l.commutable-up-toiterations}, there exist a positive integer $l$ which is a factor of  $|\det(L_f-\Id_{\RR^3})|$ and a lift $\wh G$  of $g^l$ to $\RR^3$ such that $F\circ \wh G=\wh G\circ F.$   By Lemma~\ref{l.common-fixed-point-for-lifts}, $\wh G$ admits fixed points, hence 
 $\wh G$ satisfies the assumption of Proposition~\ref{p.center-fixing},  which gives that      $ g^l=\Id_{\TT^3}.$
\endproof
\begin{Remark}~\label{r.homotopic-to-id-C1-case}
	Notice that Lemmas~\ref{l.commutable-up-toiterations} and ~\ref{l.common-fixed-point-for-lifts} are stated for $C^1$ partially hyperbolic DA diffeomorphisms. 
	By the proof of Proposition~\ref{p.center-fixing}, if $f$ is accessible, one has the same conclusion even if $f$ is $C^1$. To be precise, if $f$ is a $C^1$-partially hyperbolic diffeomorphism on $\TT^3$ whose linear part is Anosov and $f$ is accessible, then for any    $g\in\cZ^1(f)$ which is homotopic to identity, there exists a factor  $l\in\NN$ of ~$|\det(L_f-\Id_{\RR^3})|$ such that  $g^l=\Id_{\TT^3}$.
\end{Remark}

\medskip

Now, the first item in Theorem~\ref{thm.virtually-trivial-case-homotopic-to-id} is obtained. The following result completes the proof of Theorem~\ref{thm.virtually-trivial-case-homotopic-to-id}.

\begin{Proposition}~\label{p.at-most-one}
	Let $f$ be a $C^r$ ($r>1$) partially hyperbolic diffeomorphism on $\TT^3$ whose linear part $L_f$ is Anosov.  Then one has 
	$$
	\#\big\{g\in\cZ^r(f):\textrm{$g$ is homotopic to identity }\big\}\leq  |\det(L_f-\Id_{\RR^3})|.
	$$

\end{Proposition}
\proof
Let $p_1,\cdots, p_k$ be all the fixed points of $L_f$, where $k=|\det(L_f-\Id)|$.  Consider the semi-conjugacy $h:\TT^3\to\TT^3$ between $f$ and $L_f$ which is homotopic to identity.  Then all the fixed points of $f$ is contained in $\cup_{i=1}^kh^{-1}(p_i)$, and $h^{-1}(p_i)$ is an $f$-invariant center segment.
Let $I_i$ be the smallest connected segments containing all fixed points of $f$ in $h^{-1}(p_i)$. If $f$ reverses the orientation of the center foliation, then each $I_i$ is reduced to a single point. Since $\cF^c$ is orientable,  we give it an orientation. Let $I_i=[a_i,b_i]^c$ such that the direction from $a_i$ to $b_i$ gives the positive orientation. Let $\cE=\{a_i\}_{i=1}^k$. 

For any $g\in\cZ^r(f)$ which is homotopic to identity, by Corollary~\ref{co:preserve-foliations}, the center foliation $\cF^c$   is $g$-invariant and $g$ preserves the orientation of $\cF^c.$ Therefore for each $i\in\{1,\cdots, k\}$, there exists $j\in\{1,\cdots, k\}$ such that $g(I_i)=I_j$ and $g(a_i)=a_j$.
\begin{Claim-numbered}~\label{cl.identity}
	For any $g\in\cZ^r(f)$ which is homotopic to identity, if there exists some $i\in\{1,\cdots, k\}$ such that $g(a_i)=a_i$, then $g=\Id_{\TT^3}.$
	\end{Claim-numbered}
\proof[Proof of the claim]
Let $\tilde{a}_i\in\RR^3$ be lift of $a_i$. As $a_i$ is  a fixed point  for $g$, there exists a lift $G$ of   $g$ such that $G(\tilde{a}_i)=\tilde{a}_i$. Since $g\in\cZ^r(f)$ is homotopic to identity, by Proposition~\ref{p.center-fixing}, one has $G=\Id_{\RR^3},$ hence $g=\Id_{\TT^3}.$
\endproof
\begin{Claim-numbered}~\label{cl.unique-transition}
	For any $i,j\in\{1,\cdots,k\}$, there exists at most one $g\in\cZ^r(f)$ such that 
	\begin{itemize}
		\item[--] $g(a_i)=a_j$;
		\item[--] $g$ is homotopic to identity.
	\end{itemize}
\end{Claim-numbered}
\proof[Proof of the claim]
Assume  that there exist  $i,j\in \{1,\cdots, k\}$ and two diffeomorphisms $g_1,g_2\in\cZ^r(f)$ such that 
\begin{itemize}
	\item $g_1(a_i)=g_2(a_i)=a_j$;
	\item  $g_1$ and $g_2$ are homotopic to identity.
\end{itemize} 
Let $g=g_1\circ g^{-1}_2$. Then $g\in\cZ^r(f)$ is homotopic to identity and has $a_j$ as a fixed point. By Claim~\ref{cl.identity},  $g=\Id_{\TT^3},$ hence $g_1=g_2.$
\endproof
By Claims~\ref{cl.identity} and \ref{cl.unique-transition},  and the fact that $g\in\cZ^r(f)$ which is homotopic to identity must send $a_1$ to some $a_j$, one has 
$$\#\big\{g\in\cZ^r(f): \textrm{$g$ is homotopic to identity } \big\}\leq k=|\det(L_f-\Id_{\RR^3})|.$$
\endproof
\subsection{The linearization of centralizer is virtually $\ZZ^2$: rigidity case}
In this part, we discuss the case where the linear part of the centralizer  $\big\{L_g: g\in\cZ^r(f) \big\}$ is not virtually $\ZZ$.
The following theorem is the main result of this section. 


\begin{Theorem}~\label{thm.homotopic-non-trivial}
Let $f$ be a $C^r$ ($r>1$) partially hyperbolic diffeomorphism on $\TT^3$ whose linear part $L_f$ is Anosov. If there exists $g\in\cZ^r(f)$ such that  $L_g^m\notin\big\{ L_f^n\big\}_{n\in\ZZ}$ for any $m\neq 0$,  then $f$ is $C^{r-\e}$-conjugate to $L_f$ for every $\e>0$. 
\end{Theorem}
 
\begin{Remark}
	It is clear that Theorem~\ref{thm.main-thm-low-regularity} is a direct consequence of Theorems~\ref{thm.homotopic-to-id} and \ref{thm.homotopic-non-trivial}.
\end{Remark}

\medskip

The following lemma tells us that under the assumption of Theorem~\ref{thm.homotopic-non-trivial}, the linearization of the centralizer is virtually $\ZZ^2.$
\begin{Lemma}
	Let $f$ be a diffeomorphism on $\TT^3$ as in the assumption of Theorem~\ref{thm.homotopic-non-trivial}. Then the group $\big\{L_g\in GL(3,\ZZ): g\in\cZ^r(f) \big\}$ is abelian and virtually $\ZZ^2.$
\end{Lemma}
\proof
Since $f$ is partially hyperbolic, its the linear part $L_f$ has real simple spectra. By Lemma~\ref{l.rank-of-centralizer},  the group $$\big\{L_g\in GL(3,\ZZ): g\in\cZ^r(f) \big\}\subset\big\{B\in GL(3,\ZZ): L_fB=BL_f\big\}$$ is abelian and virtually $\ZZ^2$ or $\ZZ.$  Since $\cZ^r(f)$ contains elements which are virtually not homotopic to any element in $\{L_f^n\}_{n\in\ZZ}$, then $\big\{L_g\in GL(3,\ZZ): g\in\cZ^r(f) \big\}$ is virtually $\ZZ^2.$
\endproof

The following result is a corollary of Weyl Chambers picture applied to the linear part  of the centralizer $\cZ^1(f)$,   see for instance  Proposition 4.17 in \cite{HW}.
\begin{Corollary}~\label{co:another-weyl-chamber}
	Let $f$ be a $C^1$  partially hyperbolic diffeomorphism on $\TT^3$ whose linear part $L_f$ is Anosov. If there exists $g\in\cZ^1(f)$ such that  $L_g^m\notin\big\{ L_f^n\big\}_{n\in\ZZ}$ for any $m\neq 0$, then   there exists  $\tilde{g}\in\cZ^1(f)$ such that its linear part $L_{\tilde{g}}$ satisfies the following properties:  
	\begin{itemize}
		\item $L_{\tilde{g}}$ is contracting along $E^u_{L_f}$ and $E^s_{L_f}$;
		\item $L_{\tilde{g}}$ is uniformly expanding along $E^c_{L_f}$;
		\item  the splitting $E^u_{L_f}\oplus_\prec E^s_{L_f}\oplus_\prec E^c_{L_f}$ is dominated for $L_{\tilde{g}}$.
	\end{itemize}
\end{Corollary}
 A priori, one does not know if the diffeomorphism $g$ obtained in Corollary~\ref{co:another-weyl-chamber} is partially hyperbolic. To conclude, one needs further discussion.  
 
 Now, we show that the strong stable and unstable bundles are jointly integrable.
 \begin{Theorem}~\label{thm.joint-integrability}
 	Let $f$ be a $C^1$  partially hyperbolic diffeomorphism on $\TT^3$ whose linear part $L_f$ is Anosov. If there exists $g\in\cZ^1(f)$ such that  $L_g^m\notin\big\{ L_f^n\big\}_{n\in\ZZ}$ for any $m\neq 0$, then $f$ is not accessible.
 \end{Theorem}
\proof 
Recall that we always assume $L_f$  has two positive Lyapunov exponents.
By Corollary~\ref{co:another-weyl-chamber}, one can assume that $g\in\cZ^1(f)$ satisfies the following properties:
\begin{itemize}
	\item $ L_g$ is contracting along $E^u_{L_f}$ and $E^s_{L_f}$;
	\item $ L_g$ is uniformly expanding along $E^c_{L_f}$;
	\item  the splitting $E^u_{L_f}\oplus_\prec E^s_{L_f}\oplus_\prec E^c_{L_f}$ is dominated for $ L_g$.
\end{itemize}

Up to replacing $f$ and $g$ by $f^2$ and $g^2$,   one can assume that $f$ and $g$ preserve the orientation of $E^s_{L_f}, E^c_{L_f}, E^u_{L_f}.$ Let $\cF^u,\cF^s,\cF^c$ be the strong unstable, strong stable and center foliations of $f$ respectively. Let  $\cW^u,\cW^s,\cW^c$ be the strong unstable, strong stable and center foliations of $L_f$ respectively. Their corresponding center stable and center unstable foliations would be denoted as $\cF^{cs}, \cF^{cu}, \cW^{cs},\cW^{cu}$.

Consider a lift $F$ of $f$ to $\RR^3$. By Lemma~\ref{l.commutable-up-toiterations}, there exists $0<l\leq |\det(L_f-\Id_{\RR^3})|$ such that $g^l$ admits a lift to $\RR^3$ which commutes with $F$. For simplicity, we will  assume that $l=1$. Let $G$ be the lift of $g$ such that $F\circ G=G\circ F.$  By Lemma~\ref{l.common-fixed-point-for-lifts},  $F$ and $G$ have a common fixed point $p\in\RR^3.$
Let $H:\RR^3\to\RR^3$ be the continuous surjective map given by Corollary~\ref{co:linearizable-simultaneously} such that 
\begin{itemize}
	\item $H\circ F= L_f\circ H$ and $H\circ G=L_g\circ H;$
	\item  $H-\Id_{\RR^3}$ is $\ZZ^3$-periodic.
\end{itemize}
Then $H(p)=0.$

Now we show that $G$ has certain topological hyperbolicity.
\begin{Claim-numbered}~\label{cl:topo-contracting}
 	$G$ is topologically contracting along the foliation $\tilde\cF^s$, i.e. for any two points $x,y$ on a same $\tilde\cF^s$-leaf, one has $\lim_{n\rightarrow+\infty}\ud\big(G^n(x), G^n(y)\big)=0.$
	\end{Claim-numbered}
\proof[Proof of the claim]
By Theorem~\ref{thm.semi-conjugacy}, the map $H$ is injective along each leaf of $\tilde\cF^s$ and sends a leaf of $\tilde\cF^s$ to a leaf of $\tilde\cW^s$. Recall that $L_g$ is uniformly contracting along $\tilde\cW^s$ and $H\circ G=L_g\circ H$, then  one deduces that $G$ is topologically contracting along $\tilde\cF^s.$
\endproof
\begin{Claim-numbered}\label{cl:bounded-along-unstable}
	There exists $K>0$ such that for any $x\in\RR^3$ and $y\in\tilde\cF^u(x)$, one has 
	$$\limsup_{n\rightarrow+\infty}\ud_{\tilde\cF^u}\big(G^n(x),G^n(y)\big)\leq K.$$
	\end{Claim-numbered}
\proof[Proof of the claim]
 Recall that  $p\in\RR^3$ is a fixed point of $G$ and $H(p)=0$. Since the foliations $\tilde\cF^{cs}$ and $\tilde\cF^u$ have the global product structure, the space of $\tilde{\cF}^{cs}$-leaves can be identified as  $\tilde\cF^u(p)\cong\RR.$   Similarly, one can identify the space of $\tilde{\cW}^{cs}$-leaves of $L_f$ as $\tilde{\cW}^u(0)$.

 Thanks to Corollary~\ref{co:preserve-foliations}, one can consider   the action $G^{cs}$, induced by $G$, on the space of $\tilde{\cF}^{cs}$-leaves. Then we have the following commuting diagram, and $G^{cs}$ can be identified as the diffeomorphism $G:\tilde\cF^u(p)\to \tilde\cF^u(p)$.
 \begin{displaymath}
 \xymatrix{
 	\mathbb{R}^3 \ar[r]^{G} \ar[d]_{{\rm proj}_{cs}} & \mathbb{R}^3 \ar[d]^{{\rm proj}_{cs}} \\
 	\mathbb{R}^3/\tilde{\cF}^{cs} \ar[r]^{G^{cs}}  \ar[d]_{\cong}  & \mathbb{R}^3/\tilde{\cF}^{cs} \ar[d]^{\cong} \\
    \tilde{\cF}^u(p) \ar[r]^{G} & \tilde{\cF}^u(p) }
    \qquad \qquad \qquad
 \xymatrix{
  	\mathbb{R}^3/\tilde{\cF}^{cs} \ar[r]^{F} \ar[d]_{H^{cs}} & \mathbb{R}^3/\tilde{\cF}^{cs} \ar[d]^{H^{cs}} \\
  	\mathbb{R}^3/\tilde{\cW}^{cs} \ar[r]^{L_f}
  	\ar[d]_{\cong} & \mathbb{R}^3/\tilde{\cW}^{cs} \ar[d]^{\cong}  \\
    \tilde{\cW}^u(0) \ar[r]^{L_f} & \tilde{\cW}^u(0) }
 \end{displaymath}
 By Theorem~\ref{thm.semi-conjugacy},  
 the map $H$ sends a center stable leaf of $F$ to a center stable leaf of $L_f$ which induces a map $H^{cs}$ from the space of $\tilde\cF^{cs}$-leaves to the space of   $\tilde{\cW}^{cs}$-leaves.

 Since $H^{-1}(x)$ is contained in a single center leaf for any $x\in\RR^3$, then  $H^{cs}$ is a homeomorphism  from the space of $\tilde{\cF}^{cs}$-leaves to the space of $\tilde{\cW}^{cs}$-leaves.
 Combining  with the fact that $H\circ G=L_g\circ H$, one gets that the  homeomorphism $G^{cs}:\tilde{\cF}^u(p)\to\tilde{\cF}^u(p)$ is conjugate to  $L_g:\tilde{\cW}^u(0)\to \tilde{\cW}^u(0)$. 
 \begin{displaymath}
    \xymatrix{
    	\mathbb{R}^3/\tilde{\cF}^{cs}\cong\tilde{\cF}^u(p) \ar[r]^{G^{cs}} \ar[d]_{H^{cs}} & \mathbb{R}^3/\tilde{\cF}^{cs}\cong\tilde{\cF}^u(p) \ar[d]^{H^{cs}} \\
    	\mathbb{R}^3/\tilde{\cW}^{cs}\cong{\cW}^u(0) \ar[r]^{L_g}   & \mathbb{R}^3/\tilde{\cW}^{cs}\cong{\cW}^u(0)  }
 \end{displaymath}
 By the choice of $g$, the linear map $L_g$   is a contracting along $\tilde{\cW}^u(0)$, therefore $G^{cs}$ is topological contracting, so is $G:\tilde\cF^u(p)\to \tilde\cF^u(p).$
 
Let $x\in\RR^3$ and   $y\in\tilde{\cF}^u(x)$.  By the global product structure, 
 the leaves $\tilde{\cF}^{cs}(x)$ and $\tilde{\cF}^{cs}(y)$ intersect  $\tilde{\cF}^{u}(p)$ into unique points $\wh x$ and $\wh y$ respectively. Since $\ud(G^n(\wh x), G^n(\wh y))$ tends to $0$ when $n$ tends to infinity and the center stable foliation and strong unstable foliations are invariant under $G$, by Corollary~\ref{co:quantitative-product-structure}, there exists a constant $K>0$ such that 
 $$\ud_{\tilde\cF^u}\big(G^n(x), G^n(y)\big)\leq K+\ud\big(G^n(\wh x), G^n(\wh y)\big).$$  
 Letting $n$ tend  to $+\infty$, one gets the posited property.
\endproof

Assume, on the contrary, that $f$ is  accessible. Lemma~\ref{l.lift-preserve-accessibility} shows that the lift $F:\RR^3\rightarrow\RR^3$ is also accessible.
For every point $x\in\RR^3$, we choose $y\in\cF^c(x)$ such that $H(x)\neq H(y)$. Let $I^c$ denote the segment between $x$ and $y$ contained in $\cF^c(x)$.

Since $F$ is accessible, there exists a sequence of segments $I_1,I_2,\cdots,I_k$ such that for every $j=1,\cdots,k$, one has that 
\begin{itemize}
	\item $I_j$ is contained in a leaf of $\tilde{\cF}^s$ or $\tilde{\cF}^u$;
	\item The endpoints of $I_j$ are $x_{j-1}$ and $x_j$, where $x_0=x$ and $x_k=y$.
\end{itemize}
By Claim~\ref{cl:topo-contracting} and Claim~\ref{cl:bounded-along-unstable}, we have
 \begin{eqnarray}
\max_{j=1,\cdots,k}\big\{\limsup_{n\rightarrow+\infty}\ell\big(G^n(I_j)\big)\big\}\leq K, \nonumber
\end{eqnarray}
where $\ell(\cdot)$ denotes the length of a $C^1$ curve.
This implies the two endpoints of $G^n(I^c)$ are at uniformly bounded distance. 

On the other hand, since $H(x)\neq H(y)$, by the choice of $g$, one has 
$$
\lim_{n\rightarrow+\infty}\ell\big(L_g^n\circ H(I^c)\big)=\lim_{n\rightarrow+\infty}\ell\big(H\circ G^n(I^c)\big)=+\infty.
$$ 
Since $H-\Id_{\RR^3}$ is uniformly bounded on $\RR^3$, one has   
$\lim_{n\rightarrow+\infty}\ell\big(G^n(I^c)\big)=+\infty.$
This contradicts to the quasi-isometric property of the center foliation $\tilde\cF^c$ given by Theorem~\ref{thm.dynamical-coherence-unique}.
\endproof

Now, we are ready to give the proof of Theorem~\ref{thm.homotopic-non-trivial}.
\proof[Proof of Theorem~\ref{thm.homotopic-non-trivial}]
By    Theorems~\ref{thm.DA-is-ergodic} and~\ref{thm.joint-integrability}, the strong stable and unstable bundles of $f$ are jointly integrable, and $f$ is Anosov.  Let $h$ be the homeomorphism such that $h\circ f=L_f\circ h.$  

Now, let us fix some notations. Recall the partially hyperbolic splitting for $f$ is denoted as $E^s\oplus E^c\oplus E^u$ and the partially hyperbolic splitting for $L_f$ is denoted as $E^s_{L_f}\oplus E^c_{L_f}\oplus E^u_{L_f}$. Let $\cF^s,\cF^u,\cF^c$ be the strong stable, strong unstable and center foliations of $f$ respectively, and by Corollary~\ref{co:preserve-foliations}, these foliations are invariant under each element of $\cZ^r(f)$.  Let $\cW^s,\cW^u,\cW^c$ be the strong stable, strong unstable and center foliations of $L_f$ respectively. 
 
Since $E^s\oplus E^u$ is integrable, by Theorem~\ref{thm.gan-shi},  one has that 
\begin{itemize}
	\item the conjugacy $h$ is uniformly $C^1$ along the center leaves of $f$;
	\item  $h$ sends the foliation $\cF^u$ to the foliation $\cW^u.$
\end{itemize}

In the following, we will show that $h$ is uniformly $C^1$ along the leaves of $\cF^s$ and $\cF^u.$

\medskip

As the discussion in \cite{H}, since the rank of the linearization of whole $\cZ^r(f)$ action is $2$ (which is full rank),  hence  $\cZ^r(f)$ induces a maximal Cartan affine action on the torus,  and  there exist  diffeomorphisms $g, \wh g\in\cZ^r(f)$ whose   linear parts satisfy: 
\begin{itemize}
	\item $L_g$ is uniformly expanding along $E^c_{L_f}$, and $L_g$ is uniformly contracting along $E^u_{L_f}\oplus E^{s}_{L_f}$;
	\item $L_{\wh g}$ is uniformly expanding along $E^s_{L_f}\oplus E^c_{L_f}$, and $L_{\wh g}$ is uniformly  contracting along $E^u_{L_f}$;
	\item  the splitting $E^u_{L_f}\oplus E^s_{L_f}\oplus E^c_{L_f}$ is dominated for $L_g$ and $L_{\wh g}$.
\end{itemize}
 By Theorem~\ref{thm.adler-palais}, one has that $h\circ g=L_g\circ h$ and $h\circ \wh g=L_{\wh g}\circ g.$ Since  $h$ is uniformly $C^1$ along the center leaves of $f$, then  $g$ and $\wh g$ are uniformly expanding along $E^{c}_f$. 

\begin{Claim-numbered}
The diffeomorphisms $g$ and $\wh g$ are Anosov. To be precise:
\begin{itemize}
	\item $g$  is uniformly contracting along $E^u$ and $E^s$;
	\item $\wh g$ is uniformly contracting along $E^u$, and uniformly expanding along $E^s$. 
\end{itemize}	
\end{Claim-numbered}
\proof[Proof of the claim]  We only prove the case for $g$ and the case for $\wh g$ is analogously, since what we need are the conjugation through $h$  to their linear parts and the H\"older continuity of $h$ along the leaves of $\cF^u$ and $\cF^s$.

Since $g$ is topologically conjugate to $L_g$ by $h$, it satisfies the shadowing lemma. In particular, every ergodic measure of $g$ can be approximated by the atomic measures supporting on periodic orbits. Thus for proving $g$ is uniformly contracting along $E^u$ and $E^s$, it suffices to show that the Lyapunov exponents of periodic points of $g$ along $E^u$ and $ E^s$ are uniformly smaller than zero.

 By the continuity of the bundle $E^u$, for any $\e>0$, there exists $\delta>0$ such that for any   $x,y\in\TT^3$ with $\ud(x,y)<\delta$, one has 
\begin{eqnarray}~\label{eq:continuity-center-norm}
-\e\leq \log\|Dg|_{E^u(x)}\|-\log{\|Dg|_{E^u(y)}\|}\leq \e.
\end{eqnarray}
Since $g$ is conjugate to $L_g$ and $h$ sends the foliations $\cF^u,\cF^s$ to the corresponding linear foliations of $L_g$, hence $g$ is topologically contracting along the leaves of $\cF^u$ and $\cF^s$.  Let $p$ be a periodic point of period $k$.  Then $g^k|_{\cF^u(p)}:\cF^u(p)\to \cF^u(p)$ is topologically contracting and has a unique fixed point $p.$ 
Let $x\in \cF^u_{\delta}(p)$, then $g^{nk}(x)\in\cF^u_{\delta}(p)$ for any $n\in\NN.$  By Equation~\eqref{eq:continuity-center-norm}, one has 
\begin{eqnarray}~\label{eqn:contracting-at-p}
\exp\big((\chi^u(p)-\e)nk\big)   \cdot \ud(x,p)\leq \ud\big(g^{nk}(x),g^{nk}(p)\big)\leq \exp\big((\chi^u(p)+\e)nk\big) \cdot \ud(x,p),
\end{eqnarray}
where $\chi^u(p)$ is the Lyapunov exponent of
$p$ for  $g$ along the direction $E^u$. 
Since $h\circ g= L_g\circ h$, one has 
\begin{eqnarray}~\label{eqn:conjugation-at-p}
\ud\big(h\circ g^{nk}(x), h\circ g^{nk}(p)\big)=\ud\big(L_g^{nk}(h(x)), L_g^{nk}(h(p))\big)=\exp\big(\chi^u(L_g)\cdot{nk}\big)\cdot\ud\big(h(x),h(p)\big),
\end{eqnarray} 
where $\chi^u(L_g)$ is the Lyapunov exponent  of $L_g$
along $E_{L_f}^u.$ 
By Remark~\ref{r.gan-shi}, the map  $h$ is uniformly H\"older continuous along the leaves of $\cF^u$, that is, there exist $C,\alpha>0$ such that for any two points $x_1,x_2$ on the same $\cF^u$-leaf,  one has $\ud_{\cF^u}(x_1,x_2)  \leq C\cdot  \big(\ud_{\cW^u}(h(x_1),h(x_2))\big)^\alpha.$
Hence one has $ \ud(hg^{nk}(x), hg^{nk}(p))\geq C^{-1/\alpha}\cdot (\ud(g^{nk}(x), g^{nk}(p)))^{1/\alpha},$  then combining with Equations~\eqref{eqn:contracting-at-p} and~\eqref{eqn:conjugation-at-p}, for any $n\in\NN$,
 one has 
$$C^{-1/\alpha}\cdot \exp\big((\chi^u(p)-\e)nk/\alpha\big)   \cdot  \big(\ud(x,p)\big)^{1/\alpha}\leq \exp\big(\chi^u(L_g)nk\big)\cdot\ud\big(h(x),h(p)\big),$$
which implies that $\chi^u(p)-\e\leq \alpha\cdot \chi^u(L_g)$.
The arbitrariness of $\e$ and $p$ give that the Lyapunov exponents of periodic points of $g$ along $E^u$ are uniformly bounded away from zero. Analogous argument gives that the Lyapunov exponents of periodic points of $g$ along $E^s$ are uniformly bounded away from zero. Hence $g$ is also Anosov. 
\endproof



By Theorem 2.1 in \cite{PR}, which states that the codimension one stable (or unstable) foliation of a $C^r$($r>1$) codimension one Anosov diffeomorphism is $C^1$-smooth, one has that 
\begin{itemize}
	\item the unstable foliation of $f$, which is tangent to $E^c\oplus E^u$,  is $C^1$-smooth;
	\item the stable foliation of $g$, which is tangent to $E^u\oplus E^s$,   is $C^1$-smooth;
	\item  the unstable foliation of $\wh g$, which is tangent to $E^s\oplus E^c$, is $C^1$-smooth.
\end{itemize}  
As a consequence, the  foliations $\cF^{s},\cF^c, \cF^u$ are  $C^1$. 
Now, both  $g$ and $\wh g$ satisfy the assumption of Proposition~\ref{thm.generalized-gan-shi-on-equi-center-exponent}, hence $h$ is uniformly  $C^1$ along the leaves of $\cF^s$ and $\cF^u$. 
As $f,g,\wh g$ are Anosov diffeomorphisms, the leaves of $\cF^s,\cF^c,\cF^u$ are $C^r.$ By Lemma~\ref{l.bootstrap}, the map $h$ is $C^r$ along the leaves of 
$\cF^s,\cF^c,\cF^u$. Finally, Journ\'e's theorem \cite{J} shows that $h\in\diff^{r-\e}(\TT^3)$ for any $\e>0.$
\endproof

\appendix
\section{Proof of Theorem~\ref{thm.generalized-gan-shi-on-equi-center-exponent}}
~\label{s.appendix}
The aim of this section is to give the proof of Proposition~\ref{thm.generalized-gan-shi-on-equi-center-exponent} which essentially follows the argument in Section 4 of  \cite{GS}. See also \cite{G12}.
\proof[Proof of Proposition~\ref{thm.generalized-gan-shi-on-equi-center-exponent}]
Assume, on the contrary,  that there exist two periodic orbits whose Lyapunov exponents along the bundle $E^c$ are different.  For each periodic point $p$, we denote by $\chi^c(p)$ the Lyapunov exponent of $p$ along the direction $E^c.$
As $g$ is Anosov and is  uniformly expanding along the  continuous bundle  $E^c$, then there exist $0<\chi_1<\chi_2$ such that $\overline{\big\{\chi^c(p): \textrm{$p$ is a periodic point} \big\}}=[\chi_1,\chi_2].$
By shadowing lemma, for each point $x\in\TT^3$, one has 
$$\chi_1\leq \liminf_{n\rightarrow\infty}\frac{1}{n}\log\|Dg^n|_{E^c(x)}\|\leq \limsup_{n\rightarrow\infty}\frac{1}{n}\log\|Dg^n|_{E^c(x)}\|\leq \chi_2.$$
Hence for any $\e>0$, there exists an adapted metric $\|\cdot \|_{\e}$ such that   
$$\chi_1-\e\leq \log\|Dg|_{E^c(x)}\|_{\e}\leq \chi_2+\e, \textrm{ for any $x\in\TT^3$}.$$
By the continuity of the bundle $E^{c}$, there exists $\delta>0$ such that for any $z_1,z_2\in\TT^3$ with $\ud(z_1,z_2)<3\delta$, one has $$-\e<\log\|Dg|_{E^c(z_1)}\|-\log\|Dg|_{E^c(z_2)}\|<\e.$$
Now, one fixes periodic points $p,q$ such that $\chi^c(p)\leq \chi_1+\e$ and $\chi^c(q)\geq\chi_2-\e$.

For the linear Anosov map $L_g$, let us denote by $\cW^{s}, \cW^c,\cW^u,\cW^{su}$ the foliations tangent to $E^{s}_{L_g}, E^c_{L_g}, E^u_{L_g}, E^{s}_{L_g}\oplus E^u_{L_g}$ respectively.

 Since the factors of unit eigenvectors of $L_g$ are algebraic, there exists $C_1>1$ such that  for any $l>0$ large, each  strong unstable segment   of $L_g$ with length $l$ is $C_1/\sqrt{l}$ dense in $\TT^3.$ 
For any $l>0$ large, there exist $x,y\in\TT^3$ such that 
\begin{eqnarray}~\label{eqn:p-q-path}
x\in \cW^u_l(h(p)) 
\qquad \textrm{ and } \qquad
 y\in \cW^s_{C_1/\sqrt{l}}(x)\cap \cW^c_{C_1/\sqrt{l}}(h(q)).
\end{eqnarray}

By the continuity of $h$, there exists $\eta>0$ such that for any center segment for $L_g$ of length no more than $\eta$, its pre-image under $h$ has length no more than $\delta.$ We will use $\ell(I)$ denote the length of a $C^1$-curve. As $\cW^{su}$ is a linear foliation, the holonomy map given by $\cW^{su}$ is an isometry.  Now, one chooses $\cW^c$-center segments $I_x$ and $I_{h(p)}$ 
such that 
\begin{itemize}
	\item  $\ell(I_x)=\ell(I_{h(p)})=\eta;$
	\item $x$ is an endpoint of $I_x$ and $h(p)$ is an endpoint of $h(p)$;
	\item  $I_x$ is an image of $I_{h(p)}$ under the holonomy map $\cW^{su}$.
\end{itemize}
Then $J_{\wh x}=h^{-1}(I_x)$ and $J_p=h^{-1}(I_{h(p)})$ are   segments tangent to $E^c$ with length no more than $\delta$,  and $J_p$ is the image of $J_{\wh x}$ under a holonomy map of $\cF^{su},$ where  $\wh x=h^{-1}(x)$. 

By Remark~\ref{r.gan-shi}, the homeomorphism  $h$ is uniformly H\"older continuous along the leaves of $\cF^u$ and $\cF^s$,  hence there exist constants $C_2>1$ and $\theta\in(0,\frac{1}{2})$ which is only determined by $h$ such that $\wh x\in\cF_d^u(p)$ and $\wh y\in \cF^{s}_{C_2/d^{\theta}}(\wh x)\cap \cF^c_{\delta}(q)$ for $d$ large, due to Equation~\eqref{eqn:p-q-path}. As $l$ can be chosen arbitrarily large, so is $d$.

As $g$ is uniformly expanding and contracting along $E^u$ and $E^s$ respectively, then  let us denote  
 $$
 \tau=\sup_{x\in\TT^3}\|Dg^{-1}|_{E^u(x)}\|<1 
  \qquad \textrm{and }  \qquad
 \kappa=\sup_{x\in\TT^3}\log\|Dg^{-1}|_{E^s(x)}\|>1.
 $$
 Let $N_d$ be the smallest integer such that $\tau^{N_d}\cdot d\leq 1.$ Then $N_d\leq -\frac{\log d}{\log\tau}+1.$ Let $N_d^1$ be the largest integer such that 
 $$\kappa^{N_d^1}C_2/d^{\theta}\leq \delta$$
 which implies that $g^{-N_d^1}(J_{\wh x})$ is contained in the $3\delta$-neighborhood of $q.$
 Then for $d$ large enough, one has the following estimate 
 \begin{align*}
 \frac{N_d^1}{N_d}&\geq \frac{\theta\log d+\log\delta-\log C_2-1}{-\frac{\log d}{\log\tau}+1}\cdot \frac{1}{\log \kappa}
 \\
 &= \frac{\theta\log d+\log\delta-\log C_2-1}{\log d-\log\tau}\cdot \frac{-\log\tau}{\log \kappa}
 \\
 &\geq \frac{\theta}{2}\frac{-\log\tau}{\log \kappa}.
 \end{align*}
  As $g$ is uniformly expanding along $E^c$, by the choices of $J_p$ and $\delta$, one has 
  $$ \exp\big(-N_d(\chi^c(p)+\e)\big)\cdot\ell(J_p)\leq  \ell\big(g^{-N_d}(J_p)\big) $$
  and 
  $$ \ell\big(g^{-N_d}(J_x)\big)\leq \ell(J_x)\cdot \exp\big(-N_d^1(\chi^c(q)-\e)\big)\cdot \exp\big((N_d^1-N_d)\cdot(\chi_1-\e)\big).$$
 Then one has 
 \begin{align*}
 \frac{\ell(g^{-N_d}(J_x))}{\ell(g^{-N_d}(J_p))}&\leq \frac{\ell(J_x)}{\ell(J_p)} \cdot \exp(-N_d^1(\chi^c(q)-\e))\cdot \exp\big((N_d^1-N_d)(\chi_1-\e)\big)\cdot \exp(N_d(\chi^c(p)+\e))
 \\
 &\leq  \frac{\ell(J_x)}{\ell(J_p)} \cdot \exp(4N_d\e)\cdot \exp\big(N_d^1(\chi^c(p)-\chi^c(q))\big)\cdot \exp\big((-N_d^1+N_d)\cdot(\chi^c(p)-\chi_1)\big) 
 \\
 & \leq \frac{\ell(J_x)}{\ell(J_p)} \cdot \exp(8N_d\e)\cdot \exp\big(N_d^1(\chi_1-\chi_2)\big)
 \\
 &\leq  \frac{\ell(J_x)}{\ell(J_p)} \cdot \exp(8\e N_d)\cdot \exp\big(  (\chi_1-\chi_2)\cdot  \frac{\theta}{2}\frac{(-\log\tau)}{\log \kappa}\cdot N_d\big).
 \end{align*}
 
One only needs to choose $\e=\frac{1}{16}(\chi_2-\chi_1)\cdot  \frac{\theta}{2}\frac{-\log\tau}{\log \kappa},$ and one  gets that $ \frac{\ell(g^{-N_d}(J_x))}{\ell(g^{-N_d}(J_p))}$ tends to $0$ when $d$ tends to infinity. Since the holonomy map given by the foliation $\cF^u$ restricted to the unstable foliation of $f$ is uniformly $C^1$, therefore $\frac{\ell(g^{-N_d}(J_x))}{\ell(g^{-N_d}(J_p))}$ is uniformly bounded from above and below which gives the contradiction.  
This proves that all the periodic points have the same Lyapunov exponent along the bundle $E^c$. 
Applying Corollary 2.6 of \cite{GS}, one gets a periodic point whose Lyapunov exponent along $E^c$ is the same as the center Lyapunov exponent of $L_f.$
Finally, we apply Theorem 5.1 of \cite{GS}, which shows that $h$ is $C^1$ along each leaf of $\cF^c.$
\endproof

 \bibliographystyle{plain}

 \begin{tabular}{l l l}
 	\emph{\normalsize Shaobo GAN}
 	& \quad &
 	\emph{\normalsize Yi SHI}
 	\medskip\\ 	
 	\small School of Mathematical Sciences
 	&& \small School of Mathematical Sciences \\
 	\small Peking University
 	&& \small Peking University\\
 	\small Beijing, 100871, P. R. China
 	&& \small Beijing, 100871, P. R. China\\
 	\small \texttt{gansb@pku.edu.cn }
 	&& \small \texttt{shiyi@math.pku.edu.cn}\\
 	&& \small
 	\bigskip\\
 	\emph{\normalsize Disheng XU}
 	& \quad &
 	\emph{\normalsize Jinhua ZHANG}
 	\medskip\\
 	\small Department of Mathematics
 	&& \small School of Mathematical Sciences\\
 	\small Imperial College London
 	&& \small Beihang University\\
 	\small London   SW7 2AZ, UK
 	&& \small Beijing, 100191, P. R.  China\\
 	\small \texttt{dxu2@ic.ac.uk}
 	&&\small \texttt{jinhua$\_$zhang@buaa.edu.cn}\\
 	&&\small \texttt{zjh200889@gmail.com}
 	
 \end{tabular}
 
\end{document}